\documentclass[11pt]{article}
\usepackage[utf8]{inputenc}
\usepackage[T1]{fontenc}
\usepackage{amsmath,amssymb,amsthm}
\usepackage[a4paper,margin=2.5cm]{geometry}
\usepackage{hyperref}

\newtheorem{theorem}{Theorem}
\newtheorem{lemma}[theorem]{Lemma}
\newtheorem{corollary}[theorem]{Corollary}
\theoremstyle{remark}
\newtheorem{remark}[theorem]{Remark}

\newcommand{\R}{\mathbb{R}}
\newcommand{\Sgn}{\operatorname{Sgn}}
\newcommand{\sgn}{\operatorname{sgn}}
\newcommand{\Lp}[1]{|#1|_p}

\title{A necessary regularity condition for finite-time stopping\\ of the string with dry friction}
\author{Alain Haraux\\[4pt]
\small Sorbonne Universit\'e, Universit\'e Paris Cit\'e, CNRS,\\
\small Laboratoire Jacques-Louis Lions, LJLL, F-75005 Paris, France\\
\small \texttt{alain.haraux@sorbonne-universite.fr}}

\date{\today}

\begin{document}
\maketitle

\begin{abstract}
We consider the one-dimensional wave equation with dry (Coulomb) friction, $u_{tt}-u_{xx}+\Sgn(u_t)\ni 0$, with Dirichlet boundary conditions. Whether every finite-energy solution stops in finite time has been open since 1978. We show that stopping in finite time is impossible whenever the initial slope or the initial velocity is unbounded: the space $W^{1,\infty}\times L^\infty$ for $(u,u_t)$ is invariant under the flow in both time directions, while stopped states are Lipschitz. In particular the answer is negative for general finite-energy data, and the open question is reduced to Lipschitz data.
\end{abstract}

\bigskip

\qquad {\small \bf AMS classification numbers:}   35B35, 35B40, 35D30, 35L10, 35L67.
 

\section{Setting}

Let $l>0$. We consider
\begin{equation}\label{eq:main}
u_{tt}-u_{xx}+\Sgn(u_t)\ni 0 \ \text{ in } (0,\infty)\times(0,l),\qquad u(t,0)=u(t,l)=0,\qquad u(0)=u_0,\ u_t(0)=u_1,
\end{equation}
where $\Sgn(s)=\sgn(s)$ for $s\neq 0$ and $\Sgn(0)=[-1,1]$. For $(u_0,u_1)\in H^1_0(0,l)\times L^2(0,l)$ the problem generates a contraction semigroup on the energy space (Brezis~\cite{Brezis72}, maximal monotone operator). The \emph{generalized solution} $u$ satisfies
\begin{equation}\label{eq:f}
u_{tt}-u_{xx}=-f \quad\text{in } \mathcal{D}'((0,\infty)\times(0,l)),\qquad f\in L^\infty,\ |f|\le 1,
\end{equation}
where $f$ is a measurable section of $\Sgn(u_t)$.

We say that the string \emph{stops at time $T$} if $u_t(t)=0$ for all $t\ge T$. Cabannes~\cite{Cabannes78} proved stopping in finite time, with explicit stopping time, for triangular initial data: for a triangle with slopes $\alpha>0$ and $-\beta<0$ (units $l=1$, friction coefficient $1$), $T=\alpha+\beta$, in particular $T=2k$ for the symmetric tent of slope $k$. The general question had been raised in~\cite{Haraux78, HarauxThese}.

\section{Main result}

\begin{theorem}\label{thm:main}
Let $(u_0,u_1)\in H^1_0(0,l)\times L^2(0,l)$ and let $u$ be the generalized solution of~\eqref{eq:main}. For any  $T\ge0$,
\[
(u_x(T),u_t(T))\in L^\infty\times L^\infty \iff (u_0',u_1)\in L^\infty\times L^\infty .
\]
Here $u'_0$ stands for $u_x(0)$. In particular, if $(u_0',u_1)\notin L^\infty\times L^\infty$, the solution does not stop in finite time.
\end{theorem}

\begin{corollary}\label{cor:counter}
Let $\tfrac12<p<1$, $u_0(x)=(x(1-x))^p$ with $l=1$, and $u_1=0$. Then $(u_0,u_1)$ has finite energy, $u_0'$ is unbounded near $x=0$ and $x=1$, and the solution does not stop in finite time. Hence finite-time stopping fails for general finite-energy data.
\end{corollary}

\begin{remark}
The question remains open for $(u_0,u_1)\in W^{1,\infty}\times L^\infty$. The theorem says that this class is natural: it is the largest class in which stopping can occur, and it is preserved by the flow.
\end{remark}

\section{Proof}

\begin{lemma}[stopped states are Lipschitz]\label{lem:stopped}
If $u_t(t)=0$ for all $t\ge T$, then $|u_{xx}(T)|\le1$ and $|u_x(T)|\le l$.
\end{lemma}
\begin{proof}
For $t\ge T$, $u(t)=u(T)$, so~\eqref{eq:f} gives $u_{xx}(T)=f$ almost everywhere, hence $|u_{xx}(T)|\le1$. Since $u(T)$ vanishes at $x=0$ and $x=l$, $u_x(T)$ vanishes at some point, hence $|u_x(T)|\le l\sup|u_{xx}|\le l$.
\end{proof}

\begin{lemma}[uniform estimate for the linear equation with bounded source]\label{lem:linear}
Extend $u$ and $f$ to $x\in\R$ as odd, $2l$-periodic functions, and set $w_\pm=u_t\pm u_x$. For $1\le p<\infty$ let $\Lp{g}=\big(\tfrac1{2l}\int_0^{2l}|g|^p\,dx\big)^{1/p}$. Then for all $s,t\ge0$,
\begin{equation}\label{eq:est}
\Lp{w_\pm(t)}\le \Lp{w_\pm(s)}+|t-s| ,
\end{equation}
with a constant independent of $p$.
\end{lemma}
\begin{proof}
From~\eqref{eq:f}, $(\partial_t-\partial_x)w_+=-f$ and $(\partial_t+\partial_x)w_-=-f$. Thus $w_\pm$ are transported with speed $\mp1$ and perturbed by a source $-f$ with $|f|\le1$. Translations preserve $\Lp{\cdot}$ and the source integrates to at most $|t-s|$ in $\Lp{\cdot}$ , which gives~\eqref{eq:est} in both time directions. Actually the result follows immediately when $ (u(0), u_t(0)) \in (H^2\cap H^1_0)  \times H^1_0 $  and $f$ is replaced by a suitable regularization bounded with values in  $ H^1_0 $  and uniformly bounded by 1 such as $ (I-\varepsilon \partial_x^2) ^{-1} f $  with $\varepsilon>0 $  tending to $0$. Since for $p<\infty$,  
$(H^2\cap H^1_0)  \times H^1_0 $ is dense in $W^{1,p}_0\times L^p$, the estimate extends by density from regular to generalized solutions. \end{proof}

\begin{proof}[Proof of Theorem~\ref{thm:main}]
Assume $(u_0',u_1)\in L^\infty\times L^\infty$. Letting $p\to\infty$ in~\eqref{eq:est} with $s=0$, since $\Lp{g}\to\|g\|_\infty$ for $g\in L^\infty$, gives $\|w_\pm(t)\|_\infty\le\|w_\pm(0)\|_\infty+t$ for all $t\ge0$. Since $u_x=(w_+-w_-)/2$ and $u_t=(w_++w_-)/2$, $(u_x(t),u_t(t))\in L^\infty\times L^\infty$.

Conversely, assume $(u_x(T),u_t(T))\in L^\infty\times L^\infty$ for some $T\ge 0$, and set $M=\|w_\pm(T)\|_\infty$. Then $\Lp{w_\pm(T)}\le M$ for all $p$. Applying~\eqref{eq:est} with $s=T$ and $t=0$ gives $\Lp{w_\pm(0)}\le M+T$ for all $p<\infty$. A function whose normalized $L^p$ norms are bounded uniformly in $p$ belongs to $L^\infty$, with norm at most that bound. Hence $w_\pm(0)\in L^\infty$, that is $(u_0',u_1)\in L^\infty\times L^\infty$.

Finally, if the string stops at $T$, then by Lemma~\ref{lem:stopped}, $(u_x(T),u_t(T))=(u_x(T),0)\in L^\infty\times L^\infty$, so the equivalence gives $(u_0',u_1)\in L^\infty\times L^\infty$.
\end{proof}

\begin{proof}[Proof of Corollary~\ref{cor:counter}]
$u_0'(x)=p\,(x(1-x))^{p-1}(1-2x)$, which is in $L^2$ iff $2(p-1)>-1$, i.e.\ $p>\tfrac12$, and is unbounded near $0$ and $1$ iff $p<1$. Apply Theorem~\ref{thm:main}.
\end{proof}

\begin{remark}
The argument uses only the linear structure with a bounded right-hand side and holds in one space dimension. The time reversal is legitimate because the estimate~\eqref{eq:est} does not depend on the sign of the friction. For the membrane, regularity is not preserved in this way.
\end{remark}

\begin{remark}
If $(u_0',u_1)\in L^\infty\times L^\infty ,$ as a consequence of  \cite {Haraux09}, Corollaries 2.2--2.3,  we know that actually $ (u_x(t),u_t(t))$ remains bounded  uniformly for $t\ge 0$  with values in $L^\infty\times L^\infty $, this is immediate if the operator $g= \Sgn $  is replaced by its Yosida regularization  $g_\lambda$ and the result follows by letting $\lambda$ tend to $0$. But this property is only a starting point which does not allow us to conclude to finite stopping.

\end{remark}

\section{Numerical illustration}

The discrete scheme (centered finite differences in $x$, time step $0.95\,dx$, exact soft-thresholding of the friction at each step) always stops in finite time, but its stopping time depends on the mesh. For $u_0=5\,(x(1-x))^{0.6}$, $l=1$, friction coefficient $1$, the strict stopping time ($v\equiv0$) is as follows.
\begin{center}
\begin{tabular}{lccccc}
$N$ & 400 & 800 & 1600 & 3200 & 6400\\ \hline
$T$ & 31.2 & 36.7 & 43.2 & 50.7 & 59.7
\end{tabular}
\end{center}
$T$ grows like a power of $N$ (about $N^{0.23}$) with no plateau, in agreement with Theorem~\ref{thm:main}. For $p=1$ (parabola) $T=5$ on every grid, as predicted by Cabannes' Theorem~1. These experiments do not prove anything about the continuous problem.

\bigskip\noindent
{\bf Acknowledgement:} The author thanks Claude (Anthropic) for numerical experiments, discussions, and assistance in preparing the manuscript. The main idea (regularity is preserved in both time directions, while stopped states are Lipschitz) is due to the author, who has checked all statements and is solely responsible for the content. The author also considers that the discussions with Claude, and especially the numerical simulations, were fundamental for him to conceive the main idea.

\end{document}